\documentclass[12pt, reqno]{amsart}
\usepackage [latin1]{inputenc}
\usepackage{amscd}
\usepackage{epsfig}
\usepackage{amssymb}
\usepackage{amsmath}
\usepackage{amsthm}
\usepackage[T1]{fontenc}
\usepackage{ae,aecompl}

\usepackage {amsfonts} % fuer IR, IN, etc; \mathfrak (Altdeutsche Symb.) - gehoert NICHT zum Standard von 2e
\usepackage {latexsym} % noetig fuer einige Symbole u.a. \Box und \leadsto
\usepackage {exscale} % damit Summen- und Integralzeichen bei >10pt richtig angezeigt werden
\usepackage{amsmath}

\newcommand{\N} {\ensuremath{\mathbb{N}}}
\newcommand{\C} {\ensuremath{\mathbb{C}}}

\newcommand{\OO}{\mathcal{O}}

\newcommand{\dq}{\overline{\partial}}
\newcommand{\wt}[1]{\widetilde{#1}}

\DeclareMathOperator{\Sing}{Sing}

\newtheorem {satz} {Satz} [section]

\newtheorem {thm} [satz] {Theorem}

\DeclareMathOperator{\supp}{supp}

\renewcommand{\theta}{\vartheta}

\title[Hartogs' Extension Theorem] 
{A $\dq$-theoretical proof of Hartogs' Extension Theorem on Stein spaces with isolated singularities}

\author{J. Ruppenthal}

\address{Mathematisches Institut, Universit\"at Bonn, Beringstr. 1, D-53115 Bonn, Germany.}
\email{jean@math.uni-bonn.de}

\date{December 12, 2007}

\begin{document}

%~\\[-3mm]

\begin{abstract} 
Let $X$ be a connected normal Stein space of pure dimension $d\geq 2$ with isolated singularities only.
By solving a weighted $\dq$-equation with compact support on a desingularization of $X$,
we derive Hartogs' Extension Theorem on $X$ by the $\dq$-idea due to Ehrenpreis.
\end{abstract}

\maketitle

\section{Introduction}

Whereas first versions of Hartogs' Extension Theorem (e.g. on polydiscs) were usually obtained
by filling Hartogs' Figures with holomorphic discs, no such geometrical proof was known for the general Theorem
in complex number space for a long time. Proofs of the general Theorem in $\C^n$ 
usually depend on the Bochner-Martinelli-Koppelman kernel
or on the solution of the $\dq$-equation with compact support
(the famous idea due to Ehrenpreis \cite{Eh}).

Only recently, J.\ Merker and E.\ Porten were able to fill the gap by giving a geometrical proof 
of Hartogs' Extension Theorem in $\C^n$ (see \cite{MePo1}) by using a finite number of
parameterized families of holomorphic discs and Morse-theoretical tools for the global topological control 
of monodromy, but no $\dq$-theory or intergal kernels (except the Cauchy kernel).

They also extended their result to the general case of $(n-1)$-complete normal complex spaces (see \cite{MePo2}),
where no proof was known until now at all. One reason is the lack of global integral kernels or an appropriate 
$\dq$-theory for singular complex spaces. The present paper is a partial answer to the
question of Merker and Porten wether it could be possible to use some $\dq$-theoretical considerations
for reproducing their result on $(n-1)$-complex spaces. More precisely, we 
solve a weighted $\dq$-equation with compact support on a desingularization of $X$,
in order to derive the following statement by the principle of Ehrenpreis:

\begin{thm}\label{thm:hartogs}
Let $X$ be a connected normal Stein space of pure dimension $d\geq 2$ with isolated singularities only.
Furthermore, let $\Omega$ be a domain in $X$ and $K\subset \Omega$ a compact subset such that $\Omega\setminus K$ is connected.
Then each holomorphic function $f\in \mathcal{O}(\Omega\setminus K)$ has a unique holomorphic extension to the whole set $\Omega$.
\end{thm}

For a more detailed introduction to the topic, we refer to \cite{MePo1} and \cite{MePo2}.

%\newpage

\section{Proof of Theorem \ref{thm:hartogs}}

The assumption about normality implies that $X$ is reduced.
For convenience, we may assume that $X$ is embedded properly into a complex number space.\footnote{See \cite{Na1}. The result is due to R.\ Remmert}
Let $$\pi: M \rightarrow X$$ 
be a resolution of singularities,
where $M$ is a complex connected manifold of dimension $d$,
and $\pi$ is a proper holomorphic surjection.
Let $E:=\pi^{-1}(\Sing X)$ be the exceptional set of the desingularization.
Note that
\begin{eqnarray}\label{eq:bih}
\pi|_{M\setminus E}: M\setminus E \rightarrow X\setminus \Sing X
\end{eqnarray}
is a biholomorphic map. For the topic of desingularization we refer to \cite{AHL}, \cite{BiMi} and \cite{Ha}.
It follows that $M$ is a $1$-convex complex manifold, and that there exists a
strictly plurisubharmonic exhaustion function
$$\rho: M \rightarrow [-\infty,\infty)\ ,$$
such that $\rho$ takes the value $-\infty$ exactly on $E$ (see \cite{CoMi}).
We can assume that $\rho$ is real-analytic on $M\setminus E$.
Let
$$\Omega':=\pi^{-1}(\Omega),\ K':=\pi^{-1}(K),\ F:=f\circ \pi \in \mathcal{O}(\Omega'\setminus K').$$
Note that $K'$ is compact because $\pi$ is a proper map.
$\Omega\setminus K$ is a connected normal complex space.
Hence, $\Omega\setminus K\setminus\Sing X$
is still connected. So, the same is true for 
$\Omega'\setminus K'\setminus E$
because of \eqref{eq:bih}. But then $\Omega'\setminus K'$ and $\Omega'$ are connected, too.
That means that the assumptions on $\Omega$ and $K$ behave well under desingularization.
Choose $\delta>0$ such that
$$K' \subset D:=\{z\in M: \rho(z)<\delta\},$$
which is possible since $\rho$ is an exhaustion function. 
But $\rho$ is also strictly plurisubharmonic outside $E$, and 
it follows that $D$ is strongly pseudoconvex. We will use the fact that 
\begin{eqnarray}\label{eq:grauert}
\dim H^q(D, \mathcal{S}) <\infty
\end{eqnarray}
for all coherent analytic sheaves $\mathcal{S}$ and $q\geq1$.
This result goes back to Grauert who originally proved it in case $\mathcal{S}=\OO$ (see \cite{Gr1}).
Let 
$$\chi\in C^\infty_{cpt}(M)$$
be a smooth cut-off function that is identically one in a neighborhood of $K'$
and has compact support in $D\cap\Omega'$ such that
\begin{eqnarray}\label{eq:support}
\Omega'\setminus C\  \ (\mbox{with } C:=\supp \chi)
\end{eqnarray}
is connected. That is possible if we choose the neighborhhod of $K'$ small enough.\\

%\newpage
Consider
$$G:=(1-\chi) F \in C^\infty(\Omega'),$$
which is an extension of $F$ to $\Omega'$, but unfortunately not holomorphic. We have to fix it by
the idea of Ehrenpreis. So, let
$$\omega:=\dq G \in C^\infty_{0,1}(\Omega'),$$
which is in fact a $\dq$-closed $(0,1)$-form with compact support in $C \subset D\cap\Omega'$.
We are now searching a solution of the $\dq$-equation
\begin{eqnarray}\label{eq:dq}
\dq g=\omega.
\end{eqnarray}
But we only know that $D$ is a strongly pseudoconvex subset of a $1$-convex complex manifold. 
This is now the place to introduce some ideas from the $\dq$-theory on singular complex spaces.
We will use a result of Forn{\ae}ss, {\O}vrelid and Vassiliadou presented in \cite{FOV1}, Lemma 2.1.
We must verify their assumptions. So, let $$\wt{D}:=\pi(D).$$
That is a strongly pseudoconvex neighborhood of the isolated singularities in $X$.
Hence, it is a holomorphically convex subset of a Stein space by the results of Narasimhan (see \cite{Na2}),
and therefore Stein itself.

Let $\mathcal{J}$ be the sheaf of ideals of the exceptional set $E$ in $M$.
Now, the result of Forn{\ae}ss, {\O}vrelid and Vassiliadou reads as:

\begin{thm}\label{thm:fov}
Let $\mathcal{S}$ be a torsion-free coherent analytic sheaf on $D$, and $q>0$.
Then there exists a natural number $T\in\N$ such that
$$i_q: H^q(D,\mathcal{J}^T \mathcal{S}) \rightarrow H^q(D,\mathcal{S})$$
is the zero map, where $i_q$ is the map induced by the natural inclusion $\mathcal{J}^T\mathcal{S}\hookrightarrow\mathcal{S}$.
\end{thm}

This statement reflects the fact that the cohomology of $M$ is concentrated along the exceptional set $E$,
and can be killed by putting enough pressure on $E$. We will now use Theorem \ref{thm:fov} with
the choices $q=n-1$ and $\mathcal{S}=\Omega^n_D$ the canonical sheaf on $D$. So, there exists a natural number $\mu>0$
such that
\begin{eqnarray}\label{eq:fov1}
i_{n-1}: H^{n-1}(D, \mathcal{J}^{\mu} \Omega^n_D) \rightarrow H^{n-1}(D, \Omega^n_D)
\end{eqnarray}
is the zero map. We will use Serre Duality (cf. \cite{Se}) to change to the dual statement.
But, can we apply Serre-Duality to the non-compact manifold $D$? The answer is yes, because higher cohomology groups
are finite-dimensional on $D$ by the result of Grauert \eqref{eq:grauert}, and we can use Serre's criterion
(\cite{Se}, Proposition 6). So, we deduce:
\begin{eqnarray}\label{eq:fov2}
i_c: H^1_{cpt}(D, \mathcal{O}_D) \rightarrow H^1_{cpt}(D, \mathcal{J}^{-\mu} \mathcal{O}_D)
\end{eqnarray}
is the zero map, where $i_c$ is induced by the natural inclusion $\mathcal{O}_D \hookrightarrow \mathcal{J}^{-\mu}\mathcal{O}_D$.
This statement means that we can have a solution for the $\dq$-equation \eqref{eq:dq} with compact support in $D$ 
that has a pole of order $\mu$ (at most) along $E$. Let us make that precise.

\newpage
$\mathcal{J}^{-\mu}\OO_D$ is a subsheaf of the sheaf of germs of
meromorphic functions $\mathcal{M}_D$.
We will now construct a fine resolution for $\mathcal{J}^{-\mu}\OO_D$.
Let $\mathcal{C}^\infty_{0,q}$ denote the sheaf of germs of smooth $(0,q)$-forms on $D$.
We consider $\mathcal{J}^{-\mu} \mathcal{C}^\infty_{0,q}$
as subsheaves of the sheaf of germs of differential forms with measurable coefficients.
Now, we define a weighted $\dq$-operator on $\mathcal{J}^{-\mu} \mathcal{C}^\infty_{0,q}$.
Let $f \in (\mathcal{J}^{-\mu} \mathcal{C}^\infty_{0,q})_z$ for a point $z\in M$. Then $f$ can be written as $f = h^{-\mu} f_0$,
where $h\in (\OO_D)_z$ generates $\mathcal{J}_z$ and $f_0\in (\mathcal{C}^\infty_{0,q})_z$. Let
$$\dq_{-\mu} f := h^{-\mu} \dq f_0 = h^{-\mu} \dq ( h^{\mu} f).$$
We obtain the sequence
$$0 \rightarrow \mathcal{J}^{-\mu}\OO_D \hookrightarrow
\mathcal{J}^{-\mu} \mathcal{C}^\infty_{0,0} \xrightarrow{\ \dq_{-\mu}\ }
\mathcal{J}^{-\mu} \mathcal{C}^\infty_{0,1} \xrightarrow{\ \dq_{-\mu}\ }
\cdots \xrightarrow{\ \dq_{-\mu}\ }
\mathcal{J}^{-\mu} \mathcal{C}^\infty_{0,d} \rightarrow 0 ,$$
which is exact by the Grothendieck-Dolbeault Lemma and well-known regularity results.
It is a fine resolution of $\mathcal{J}^{-\mu} \OO_D$ since the $\mathcal{J}^{-\mu} \mathcal{C}^\infty_{0,q}$
are closed under multiplication by smooth cut-off functions.
Therefore, the abstract Theorem of de Rham implies that:
$$H^q(D,\mathcal{J}^{-\mu}\OO_D) \cong \frac{\mbox{ker }
(\dq_{-\mu}: \mathcal{J}^{-\mu}\mathcal{C}^\infty_{0,q}(D) \rightarrow \mathcal{J}^{-\mu}\mathcal{C}^\infty_{0,q+1}(D))}
{\mbox{Im }
(\dq_{-\mu}: \mathcal{J}^{-\mu}\mathcal{C}^\infty_{0,q-1}(D) \rightarrow \mathcal{J}^{-\mu}\mathcal{C}^\infty_{0,q}(D))},$$
and we have the analogous statement for forms and cohomology with compact support.
Recall that
$$\omega\in C^\infty_{0,1}(D)$$
is $\dq$-closed with compact support in $D$. By the natural inclusion, we have that
$$\omega\in \mathcal{J}^{-\mu} C^\infty_{0,1}(D),$$
too, and it is in fact $\dq_{-\mu}$-closed. But then \eqref{eq:fov2} tells us that there exists a solution
$g\in \mathcal{J}^{-\mu} \mathcal{C}^\infty(D)$ such that
$$\dq_{-\mu} g = \omega,$$
and $g$ has compact support in $D$. So, $g\in C^\infty(M\setminus E)$ with support
in $D$, and 
$$\dq g =\omega = \dq G =\dq \big((1-\chi)F\big)\ \ \mbox{ on } \Omega'\setminus E.$$
Recall that 
$$\supp\omega \subset C=\supp \chi,$$
and that $\Omega'\setminus C$ is connected \eqref{eq:support}. But then $M\setminus C$ 
and $M\setminus C \setminus E$ are connected,
and $g$ is a holomorphic function on $M\setminus C\setminus E$ with support in $D$. Hence
$$\supp g \subset C.$$
So,
$$\wt{F} := (1-\chi) F - g \in \OO(\Omega'\setminus E)$$
equals $F$ on $\Omega' \setminus E \setminus C$.

%\newpage
But then
$$\wt{f} := \big((1-\chi) F - g\big) \circ \pi^{-1} \ \in \mathcal{O}(\Omega\setminus \Sing X)$$
equals $f$ on some open set, and it has an extension to the whole domain $\Omega$
by Riemann's Extension Theorem for normal spaces (see \cite{GrRe}, for example).
The extension $\wt{f}$ is unique because $\Omega\setminus K$ is connected and
$X$ is globally and locally irreducible (see again \cite{GrRe}). That completes the proof
of Theorem \ref{thm:hartogs}. With a little more effort, one can remove the assumption
that $X$ should contain only finitely many singularities, because $K$ has a neighborhhod in $\Omega$
that contains only a finite number of singular points.

\vspace{5mm}
{\bf Acknowledgments}

\vspace{2mm}
This work was done while the author was visiting the University of Michigan at Ann Arbor,
supported by a fellowship within the Postdoc-Programme of the German Academic Exchange Service (DAAD).
The author would like to thank the SCV group at the University of Michigan for its hospitality.

%\newpage

\vspace{10mm}

\end{document}